\documentstyle[12pt]{amsart}

\numberwithin{equation}{section}

\newcommand{\bthm}[2]{\vskip 8pt\bf #1\hskip 2pt\bf#2\it \hskip 8pt}
\newcommand{\ethm}{\vskip 8pt\rm}
\newcommand{\s}{\,\,\,\,}
\newcommand{\proof}{\it Proof.  \hskip 5pt \rm}
\def\dint{\displaystyle{\int}}

\begin{document}

\title[Solutions for Toda systems on Riemann surfaces]
{Solutions for Toda systems on Riemann surfaces}

\author{Jiayu Li, Yuxiang Li}

\thanks{The research was supported by  NSFC }

\address{Math. Group, The abdus salam ICTP\\ Trieste 34100,
   Italy\\
   and Academy of Mathematics and Systems Sciences\\ Chinese Academy of
Sciences\\ Beijing 100080, P. R. of China. } \email{jyli@@ictp.it}

\address{Math. Group, The abdus salam ICTP\\ Trieste 34100,
   Italy}
\email{liy@@ictp.trieste.it}


\begin{abstract}
In this paper, we study the solutions of Toda systems on Riemann
surface in the critical case, we prove a sufficient condition for
the existence of solutions of Toda systems.
\end{abstract}

\maketitle

{Dedicated to Professor Ding Weiyue on the occasion of his 60's
birthday.}

\section{Introduction}
Let $(\Sigma, g)$ be a compact Riemann surface with unit area $1$.
Ding-Jost-Li-Wang \cite{DJLW} studied the differential equation
$\Delta u = 8\pi - 8\pi he^u$ on $(\Sigma, g)$, it is the so
called Kazdan-Warner problem \cite{KW} related to the Abelian
Chern-Simons model (see \cite{CY1}, \cite{ CY2}, \cite{CD},
\cite{Clin}, \cite{CL}, \cite{CLMP1}, \cite{CLMP2},  \cite{DJLW2},
\cite{DJLW3}, \cite{DJLPW}, \cite{ST}, \cite{HKP}, \cite{JaW},
\cite{Ta}, \cite{NT}, etc). They pursued a variational approach to
the problem. They tried to minimize the functional
\begin{equation}J(u)=\frac{1}{2}\int_\Sigma|\nabla u|^2dV_g+8\pi \int_\Sigma
udV_g-8\pi \log \int_\Sigma he^udV_g \geq C, ~\text{in }
H^{1,2}(\Sigma),\end{equation} for some constant $C>0$.  Because
it is the critical case of the Moser-Trudinger inequality (1.1),
the analysis is subtle.

Let $K$ denote the Cartan matrix for $SU(N+1)$, i.e.,
$$K=(a_{ij})=
             \left(\begin{array}{rrrrrr}
                     2 & -1 & 0 & \cdots & \cdots & 0\cr
                     -1& 2 &-1& 0 & \cdot & 0\cr
                     0&-1&2&-1& \cdots & 0\cr
                     \cdots&\cdots&\cdots&\cdots&\cdots&\cdots\cr
                     0 &\cdots&\cdots& -1&2&-1\cr
                     0&\cdots&\cdots & 0 &-1&2\cr
             \end{array}\right).$$
In this paper, we consider the Toda systems on $(\Sigma, g)$ which
is related to the non Abelian Chern-Simons model \cite{NT1}:
$$-\Delta u_i = M_i(\frac {\exp (\sum_{j=1}^N a_{ij} u_j)}
{\int_\Sigma \exp (\sum_{j=1}^N a_{ij} u_j)}-1), \text { for }
1\le i \le N.$$ If $M_i<4\pi$, Jost-Wang \cite{JW} proved the
existence of solutions. In the case that $\Sigma$ is a torus,
$N=2$, $\max\{M_1,M_2\}>4\pi$ and $\min\{M_1,M_2\}\not=4\pi$,
Marcello-Margherita \cite{MM} proved the existence of the
solution.

They studied the problem by considering the functional for
$u_1,\cdots ,u_N\in H^{1,2}(\Sigma)$,

\begin{eqnarray}\Phi_{(M_1,\cdots ,M_N)}(u_1,\cdots ,u_N)&=& \frac 12 \sum_{i,j=1}^N\int_\Sigma
a_{ij}(\nabla u_i \nabla u_j+2M_iu_j)dV_g \nonumber\\
&&-\sum_{i=1}^N
 M_i\log\int_\Sigma \exp (\sum_{j=1}^N a_{ij}u_j)dV_g.
 \end{eqnarray} Jost-Wang \cite{JW} proved that the functional has a lower bound
if and only if
$$M_i\le 4\pi, \quad \text { for }i=1,2,\cdots, N.\label{cond2}$$
Marcello-Margherita \cite{MM} obtained a non-minimizing critical
point of the functional motivated by an earlier paper of
Struwe-Tarantello \cite{ST}. The idea was later also used by Djadli and Malchiodi
\cite{DJM} to study the existence of conformal metrics with
constant $Q$-curvature. It is clear that $M_i=4\pi$ is the
critical case of the functional. Whether it admits minimizer is
subtle. In this paper we study this problem. For simplicity, we
consider only the case that $N=2$, the general case need only more
calculations. In our case the functional is
$$
\Phi(u_1,u_2)= \frac 12 \sum_{i,j=1}^2\int_\Sigma a_{ij}(\nabla
u_i \nabla u_j+8\pi u_j)dV_g -\sum_{i=1}^2
 4\pi\log\int_\Sigma \exp (\sum_{j=1}^2 a_{ij}u_j)dV_g,
 $$
the Toda systems is
$$-\Delta u_i = 4\pi(\frac {\exp (\sum_{j=1}^2 a_{ij} u_j)}
{\int_\Sigma \exp (\sum_{j=1}^2 a_{ij} u_j)}-1), \text { for }
1\le i \le 2,\eqno(1.3)$$ where $a_{11}=a_{22}=2$ and
$a_{12}=a_{21}=-1$.

Our main result is as follows:

\bthm{Main Theorem}{}Let $\Sigma$  be a compact Riemann surface
with area $1$. If the Gauss curvature $K$ of $\Sigma$ satisfies
that
$$\max_{p\in\Sigma}K(p)<2\pi ,\eqno (1.4)$$
then $\Phi(u_1,u_2)$ has a minimizer.\ethm

We consider the sequence of minimizers $u^\epsilon=(u^\epsilon_1,
u^\epsilon_2)$ of $\Phi_{(4\pi-\epsilon,4\pi-\epsilon)}$ for small
$\epsilon>0$. Then $u^\epsilon$ satisfies a Toda type system. If
$u^\epsilon$ converges to $u^0=(u_1^0,u_2^0)$ in
$H_2:=H^{1,2}(\Sigma)\times H^{1,2}(\Sigma)$, then it is clear
that $\Phi(u^0)=\inf_{u\in H_2}\Phi(u)$, i.e., $u^0$ is a
minimizer of $\Phi$. If $u^\epsilon$ does not converge in $H_2$,
in this case, we say that $u^\epsilon$ blows up. Then there are
two cases happened according to Jost-Wang's result. For each case,
we derive a delicate lower bound of $\Phi$ which is one of the
main points in this paper. We apply capacity to calculate the
lower bound, so that we need not know details in the neck. Such a
trick has been used by the second author of this paper in
\cite{L}, \cite{L2} to prove the existence of extremal functions
for the classical Moser-Trudinger inequality on a compact
manifold. Another main point of this paper is the delicate
constructions of blowing up sequences $\phi^\epsilon$ in both
cases, so that $\Phi(\phi^\epsilon)$ are strictly less than the
lower bound derived before, and consequently we get a
contradiction to the assumption that $u^\epsilon$ blows up, which
proves our main theorem.

\section{Review of known results}
For any $u=(u_1,u_2)\in H^{1,2}(\Sigma)\times H^{1,2}(\Sigma)$, we
set
\begin{eqnarray*}
 \Phi_{\epsilon}(u)&=&\frac{1}{3}\int_{\Sigma}\left(|\nabla u_1|^2+|\nabla u_2|^2
   +\nabla u_1\nabla u_2
+3(4\pi-\epsilon) u_1+3(4\pi-\epsilon) u_2\right)dV_g\\
&&-(4\pi-\epsilon)\log{\int_{\Sigma}e^{u_1}}dV_g -(4\pi-\epsilon)
\log{\int_{\Sigma}e^{u_2}}dV_g.
\end{eqnarray*} It is not
difficult to check that
$$\Phi_{(4\pi-\epsilon,4\pi-\epsilon)}(v)=\Phi_\epsilon(u),$$ if we
set $v_1=\frac{2u_1+u_2}{3}$ and $v_2=\frac{u_1+2u_2}{3}$.

By Jost-Wang's result (\cite{JW} Corollary 4.6), one sees that
$\Phi_{\epsilon}$ has a minimizer $u^\epsilon$ of the functional
$\Phi_{\epsilon}(u)$ , i.e. we can find $u^\epsilon\in
H^{1,2}(\Sigma)\times H^{1,2}(\Sigma)$ such that
$$\Phi_\epsilon(u^\epsilon)=\inf\Phi_\epsilon(u).$$
Without loss of generality, we may assume that
$$\dint e^{u_1^\epsilon}dV_g=\dint e^{u_2^\epsilon}dV_g=1.$$
Then, we have the following equations:
$$\left\{\begin{array}{l}
            -\Delta u_1^\epsilon=(8\pi-2\epsilon)e^{u_1^\epsilon}-(4\pi-
              \epsilon)e^{u_2^\epsilon}-(4\pi-\epsilon)\\[1.7ex]
            -\Delta u_2^\epsilon=(8\pi-2\epsilon)e^{u_2^\epsilon}-(4\pi-
              \epsilon)e^{u_1^\epsilon}-(4\pi-\epsilon)
          \end{array}\right.$$

For $i=1,2$, let
$$S_i=\{x\in\Sigma: \hbox{there is a sequence }y^\epsilon\rightarrow x
\hbox{ s.t. } u_i^\epsilon(y^\epsilon)\rightarrow+\infty\}.$$

Jost-Wang \cite{JW} (section 5) proved that, there will be two
possibilities:

\noindent{\bf case 1:} $S_1=\{p_1\}$, and $S_2=\{p_2\}$, where
$p_1$ $p_2$ are two different points in $\Sigma$.

In this case, we set, for $i=1,2$,
$$m_i^\epsilon=u_i^\epsilon(x_i^\epsilon)=\max u_i^\epsilon,\s (r_i^\epsilon)^2
=e^{-m_i^\epsilon},\s \bar{u}_i^\epsilon=\dint_\Sigma u_i^\epsilon
dV_g.$$ Let $(\Omega_i, x=(x^1,x^2))$ be an isothermal coordinate
system around $p_i$ ($i=1,2$), and we assume the metric
$$g|_{\Omega_i}=e^{\varphi_i}((dx^1)^2+(dx^2)^2)$$
with $\varphi_i(0)=0$, $i=1,2$.

We set, for $i=1,2$, $\Omega_i^\epsilon=\{x\in {\mathbb R}^2:
x_i^\epsilon+r_i^\epsilon x\in\Omega_i\}$, which expands to the
whole ${\mathbb R}^2$. In $\Omega_1^\epsilon$, we have the
equations: \begin{eqnarray*}-\Delta_0
(u_1^\epsilon(x_1^\epsilon+r_1^\epsilon x)-m_1^\epsilon)
&=&e^{-\varphi_1(x_1^\epsilon+r_1^\epsilon x)}\left(
(8\pi-2\epsilon)e^{u_1^\epsilon(x_1^\epsilon+r_1^\epsilon
x)-m_1^\epsilon}\right.\\&&\left.- (r_1^\epsilon)^2(4\pi
-\epsilon)e^{u_2^\epsilon(x_1^\epsilon+r_1^\epsilon
x)}-(r_1^\epsilon)^2(4\pi- \epsilon)\right),\end{eqnarray*} where
$-\Delta_0=\frac{\partial^2}{\partial^2x^1}+\frac{\partial^2}{\partial^2x^2}$.
Since  $u_2^\epsilon$ are bounded from above in
$\Omega_1^\epsilon$, it follows from the Harnack inequality and
the elliptic estimates that $u_1^\epsilon$ converges in $C_{\rm
loc}^k({\mathbb R}^2)$ for any $k$ to the function $w$ which
satisfies the equation
$$\left\{\begin{array}{l}
             -\Delta w=8\pi e^w,\s\s\forall x\in{\mathbb{R}}^2\\
             w(x)\leq w(0)=0,\s and\s \dint_{{\mathbb{R}}^2}e^wdx\leq 1.
         \end{array}\right.$$
Hence, by the result in \cite{CL}, we know that
$$w=-2\log(1+\pi|x|^2).$$
In the same way,
$u_2^\epsilon(x_2^\epsilon+r_2^\epsilon x)-m_2^\epsilon$ converges to $w$.

We set $\bar{u}_i^\epsilon=\int_\Sigma u_i^\epsilon dV_g$, we have the following proposition
(see Lemma 5.6, and the proof of Theorem 3.1 in \cite{JW}).

\bthm{Proposition}{2.1} We have $\bar u^\epsilon_j\rightarrow-\infty$ for $j=1,2$.
Furthermore, for any $q\in (1,2)$, we have
$$u^\epsilon_j- \bar u^\epsilon_j \text{ converges to } G_j \text{ in } H^{1,q}(\Sigma),$$
where $G_1$ and $G_2$ satisfy
$$\left\{ \begin{array}{crl}
     -\Delta G_1& =& 8 \pi \delta_{p_1} -4\pi \delta_{p_2}-4\pi,\\
     -\Delta G_2& =& 8 \pi \delta_{p_2} -4\pi \delta_{p_1}-4\pi,\\
     \int_\Sigma G_jdV_g& =& 0, \quad \text{for } j=1,2
  \end{array}\right. $$
where $\delta_y$ is the Dirac distribution. Moreover,
$$u^\epsilon_j- \bar u^\epsilon_j \text{ converges to } G_j \text{ in } C^2_{loc}
(\Sigma{\backslash}
\{p_1,p_2\}).$$\ethm

{\bf Remark 2.1:} It is easy to see that, in $\Omega_1$,
$$G_1=-4\log{r}+A_1(p_1)+f_1,\s and\s G_2=2\log{r}+A_2(p_1)+g_1\eqno(2.1)$$
where $r^2=x_1^2+x_2^2$, $A_i(p_1)$ ($i=1,2$) are constants, and
$f_1$, $g_1$ are smooth functions which are zero at $0$.
Similarly, in $\Omega_2$, we can write
$$G_1=2\log{r}+A_1(p_2)+f_2,\s\s and\s G_2=-4\log{r}+A_2(p_2)+g_2.\eqno(2.2)$$
where $A_i(p_1)$ ($i=1,2$) are constants, and $f_2$, $g_2$ are
smooth functions which are zero at $0$. \vspace{2ex}

\noindent {\bf Case 2:} $S_1=\{p\}$, and $S_2=\emptyset$.

In this case, $u_2^\epsilon$ are bounded from above. Let
$(\Omega;x)$ be an isothermal coordinate system around $p$,
similar to the case 1, we have
$$u_1^\epsilon(x_1^\epsilon+r_1^\epsilon x)-m_1^\epsilon\rightarrow-2\log{(1+\pi|x|^2)}.$$
We also have the following proposition (c.f. \cite{JW}):

\bthm{Proposition}{2.2} Let $\bar u^\epsilon_1$ be the average of
$u^\epsilon_1$. We have $\bar u^\epsilon_1\rightarrow-\infty$.
Furthermore, for any $q\in (1,2)$, we have
$$u^\epsilon_1- \bar u^\epsilon_1 \text{ converges to } G_1 \text{ in } H^{1,q}(\Sigma),$$
and
$$u^\epsilon_2 \text{ converges to } G_2 \text{ in } H^{1,q}(\Sigma),$$
where $G_1$ and $G_2$ satisfy
$$\left\{ \begin{array}{crl}
     -\Delta G_1& =& 8 \pi \delta_{p} -4\pi e^{G_2}-4\pi,\\
     -\Delta G_2& =& 8 \pi e^{G_2} -4\pi \delta_p-4\pi,\\
     \int_\Sigma G_1dV_g& =& 0, \s\int_\Sigma e^{G_2}dV_g=1,\s\sup\limits_{x\in\Sigma}G_2<+\infty
  \end{array}\right.\eqno(2.3) $$
where $\delta_y$ is the Dirac distribution. Moreover,
$$u^\epsilon_1- \bar u^\epsilon_1 \text{ converges to } G_1,\text{ and }
u^\epsilon_2\text{ converges to } G_2 \text{ in } C^2_{loc}
(\Sigma{\backslash}
\{p\}).$$\ethm

Since $G_2$ is bounded from above, we can deduce from the equation
(2.3) that $G_2=2\log{r}+h$ in $\Omega$, where $h\in
H^{2,q}_{loc}(\Omega)$ for any $q>0$. Then  $e^{G_2}=r^2e^h\in
C^1_{loc}(\Omega)$, and then $\Delta_0 h\in C^1_{loc}(\Omega)$.
Therefore, by the standard elliptic estimates, $G_2-2\log{r}$ is
smooth in $\Omega$. So, we can write
$$G_1=-4\log{r}+A_1(p)+f,\s\s and\s G_2=2\log{r}+A_2(p)+g\eqno(2.4)$$
where $r^2=x_1^2+x_2^2$, $A_i(p)$ ($i=1,2$) are constants and $f$,
$g$ are smooth functions which are zero at $0$.

\section{The lower bound for case 1}

We assume that $\Omega_1\cap\Omega_2=\emptyset$, and
$B_r(p_1)\subset\Omega_1$. We set
$v_2^\epsilon=\frac{1}{3}(2u_2^\epsilon+u_1^\epsilon)-\frac{1}{3}
(2\bar{u}_2^\epsilon+\bar{u}_1^\epsilon)$. Then, in $B_r(p_1)$, we
have
$$\left\{\begin{array}{l}
             -\Delta v_2^\epsilon
               =(4\pi-\epsilon)e^{u_2^\epsilon}-(4\pi-\epsilon)\in L^\infty(B_r(p_1))\\[1.7ex]
             v_2^\epsilon|_{\partial B_r(p_1)}\rightarrow
             \frac{1}{3}(2G_2+G_1)
         \end{array}\right.$$
So $\|v_2^\epsilon\|_{C^1}\leq C$, where $C$ is a constant
depending only on $r$.

By a direct calculation, one gets,
\begin{eqnarray*}\frac{1}{3}\dint_{B_\delta(x_1^\epsilon)}(|\nabla
u_1|^2+|\nabla u_2|^2 +\nabla u_1\nabla
u_2)dV_g&=&\frac{1}{4}\dint_{B_\delta(x_k)}( |\nabla
u_1^\epsilon|^2+3|\nabla v_2^\epsilon|^2)dV_g\\
&= &\frac{1}{4}\dint_{B_\delta(x_1^\epsilon)}|\nabla
u_1^\epsilon|^2dV_g+O(\delta^2).\end{eqnarray*}
Recall that
$u_1(x_1^\epsilon+r_1^\epsilon x)-m_1^\epsilon\rightarrow w$ in
$C^k(B_L(0))$, for any $k$, we have
$$\frac{1}{4}\dint_{B_\delta(x_1^\epsilon)}|\nabla u_1^\epsilon|^2dV_g=
\frac{1}{4}\dint_{B_L}|\nabla w|^2dx+\frac{1}{4}\dint_{B_\delta(x_1^\epsilon)\setminus B_{Lr_1^\epsilon}(x_1^\epsilon)}
|\nabla u_1^\epsilon|^2dV_g+o(1)+O(\delta^2).$$
Let
$$a_1^\epsilon=\inf\limits_{\partial B_{Lr_1^\epsilon}(x_1^\epsilon)}u_1^\epsilon
,\s b_1^\epsilon=\sup\limits_{\partial
B_\delta(x_1^\epsilon)}u_1^\epsilon.$$ We set
$a_1^\epsilon-b_1^\epsilon=m_1^\epsilon-\bar{u}_1^\epsilon+d_1^\epsilon$.
It is clear that, for fixed $L$ and $\delta$,
$$d_1^\epsilon\rightarrow w(L)-\sup_{\partial B_\delta(p_1)}G_1\s as\s\epsilon\rightarrow 0.$$

Let $f_1^\epsilon=\max\{\min\{u_1^\epsilon,a_1^\epsilon\},b_1^\epsilon\}$. We get
$$\begin{array}{ll}
     \dint_{B_\delta(x_1^\epsilon)\setminus B_{Lr_1^\epsilon}(x_1^\epsilon)}|\nabla u_1^\epsilon|^2dV_g
         &\geq \dint_{B_\delta(x_1^\epsilon)\setminus B_{Lr_1^\epsilon}(x_1^\epsilon)}|\nabla
         f_1^\epsilon|^2dV_g\\[1.7ex]
         &= \dint_{B_\delta(x_1^\epsilon)\setminus B_{Lr_1^\epsilon}(x_1^\epsilon)}|\nabla_0 f_1^\epsilon|^2dx\\[1.7ex]
         &\geq\inf\limits_{\Psi|_{\partial B_{Lr_1^\epsilon}(0)}=a_1^\epsilon,
             \Psi|_{\partial B_{\delta}(0)}=b_1^\epsilon}
                \dint_{B_{Lr_1^\epsilon}(0)\setminus B_\delta(0)}|\nabla_0\Psi|^2dx
  \end{array}$$
Here, $|\nabla_0 g|^2=|\frac{\partial g}{\partial x_1}|^2+
|\frac{\partial g}{\partial x_2}|^2$.  It is well-known that
$\inf\limits_{\Psi|_{\partial B_{Lr_1^\epsilon}}=a_1^\epsilon,
             \Psi|_{\partial B_\delta}=b_1^\epsilon}
                \dint_{B_\delta\setminus B_{Lr_1^\epsilon}}|\nabla_0\Psi|^2dx$ is uniquely attained by
the function $\phi$ which satisfies the equation
$$\left\{\begin{array}{l}
             -\Delta_0 \phi=0\\
             \phi|_{\partial B_{Lr_1^\epsilon}}=a_1^\epsilon,\phi|_{\partial B_\delta}=b_1^\epsilon.
         \end{array}\right.$$
Hence,
$$\phi=\frac{a_1^\epsilon-b_1^\epsilon}{-\log{Lr_1^\epsilon}+\log{\delta}}\log{r}+
\frac{a_1^\epsilon\log\delta-b_1^\epsilon\log{Lr_k}}{-\log{Lr_1^\epsilon}+\log{\delta}},$$
and then
$$\dint_{B_\delta(0)\setminus B_{Lr_1^\epsilon}(0)}|\nabla_0 \phi|^2dx=\frac{4\pi(a_1^\epsilon-
b_1^\epsilon)^2}
{-\log(Lr_1^\epsilon)^2+\log{\delta^2}}.$$
Therefore, we have
$$\dint_{B_\delta(x_1^\epsilon)\setminus B_{Lr_1^\epsilon}(x_1^\epsilon)}|\nabla u_1^\epsilon|^2dV_g
\geq\frac{4\pi
            (m_1^\epsilon-\bar{u}_1^\epsilon+d_1^\epsilon)^2}{-\log{L^2}-
            \log{(r_1^\epsilon)^2}+\log{\delta^2}}.$$

Recall that $-\log{(r_1^\epsilon)^2}=m_1^\epsilon$, we get
$$\begin{array}{ll}
        \dint_{B_\delta(x_1^\epsilon)\setminus B_{Lr_1^\epsilon}(x_1^\epsilon)}|\nabla u_1^\epsilon|^2dV_g&\geq
          4\pi\frac{(m_1^\epsilon-\bar{u}_1^\epsilon+d_1^\epsilon)^2}{m_1^\epsilon}
          (1-\frac{\log{L^2}-\log{\delta^2}}{m_1^\epsilon})^{-1}\\[1.7ex]
         &\geq 4\pi\frac{(m_1^\epsilon-\bar{u}_1^\epsilon+d_1^\epsilon)^2}{m_1^\epsilon}
          (1+\frac{\log{L^2}-{\log{\delta^2}}}{m_1^\epsilon}+
          \frac{A}{(m_1^\epsilon)^2})\\[1.7ex]
         &\geq 4\pi\frac{(m_1^\epsilon-\bar{u}_1^\epsilon)^2}{m_1^\epsilon}+
           8\pi d_1^\epsilon(1-\frac{\bar{u}_1^\epsilon}{m_1^\epsilon})\\
           &+4\pi
           (1-\frac{\bar{u}_1^\epsilon}{m_1^\epsilon})^2(\log{L^2}-\log{\delta^2})+\frac{A'
           \bar{u}_1^\epsilon}
           {(m_1^\epsilon)^2},
  \end{array}$$
where $A$ and $A'$ are constants which depend only on $\delta$ and
$L$.

Then, we have
$$\begin{array}{l}\frac{1}{3}\dint_{B_\delta(x_1^\epsilon)}(|\nabla u_1|^2+|\nabla u_2|^2
+\nabla u_1\nabla u_2)dV_g\\[1.7ex]
 \s \geq \frac{1}{4}\dint_{B_L}|\nabla w|^2dx+
              \pi\frac{(m_1^\epsilon-\bar{u}_1^\epsilon)^2}{m_1^\epsilon}+
           2\pi d_1^\epsilon(1-\frac{\bar{u}_1^\epsilon}{m_1^\epsilon})+\pi
           (1-\frac{\bar{u}_1^\epsilon}{m_1^\epsilon})^2(\log{L^2}-\log{\delta^2})\\[1.7ex]
 \s\s+\frac{A'(\delta,L)
           \bar{u}_1^\epsilon}
           {(m_1^\epsilon)^2}+o(1)+O(\delta^2).
\end{array}$$

Similarly, we have
\begin{eqnarray*}
\lefteqn{\frac{1}{3}\dint_{B_\delta(x_2^\epsilon)}(|\nabla
u_1|^2+|\nabla u_2|^2
+\nabla u_1\nabla u_2)dV_g}\\
& \geq &\frac{1}{4}\dint_{B_L}|\nabla w|^2dx+
              \pi\frac{(m_2^\epsilon-\bar{u}_2^\epsilon)^2}{m_2^\epsilon}+
           2\pi d_1^\epsilon(1-\frac{\bar{u}_2^\epsilon}{m_2^\epsilon})\\
           &&+\pi
           (1-\frac{\bar{u}_2^\epsilon}{m_1^\epsilon})^2(\log{L^2}-\log{\delta^2})+\frac{A'
           \bar{u}_2^\epsilon}
           {(m_2^\epsilon)^2}+o(1)+O(\delta^2).
\end{eqnarray*}
It concludes that
$$\begin{array}{ll}
    \frac{1}{3}\dint_{B_\delta(x_1^\epsilon)\cup B_\delta(x_2^\epsilon)
       }(|\nabla u_1|^2+|\nabla u_2|^2
         +\nabla u_1\nabla u_2)dV_g+(4\pi-\epsilon)
          \bar{u}_1^\epsilon+(4\pi-\epsilon)\bar{u}_2^\epsilon\\[1.7ex]
    \s\geq\frac{1}{3}\dint_{B_\delta(x_1^\epsilon)\cup B_\delta(x_2^\epsilon)
       }(|\nabla u_1|^2+|\nabla u_2|^2
         +\nabla u_1\nabla u_2)dV_g+4\pi
          \bar{u}_1^\epsilon+4\pi\bar{u}_2^\epsilon\\[1.7ex]
    \s\geq \frac{1}{2}\dint_{B_L}|\nabla w|^2dx+
       \sum_{i=1,2}\left(\pi\frac{(m_i^\epsilon+\bar{u}_i^\epsilon)^2}{m_i^\epsilon}
         +2\pi d_i^\epsilon(1-\frac{\bar{u}_i^\epsilon}{m_i^\epsilon})+
          \pi(1-\frac{\bar{u}_2^\epsilon}{m_1^\epsilon})^2(\log{L^2}-
          \log{\delta^2})\right)\\[1.7ex]
    \s\s +\sum\limits_{i=1,2}\frac{A'\bar{u}_i^\epsilon}
           {(m_i^\epsilon)^2}+o(1)+O(\delta^2)\\[1.7ex]
    \s\geq \frac{1}{2}\dint_{B_L}|\nabla w|^2dx+
       \sum_{i=1,2}\left(\pi m_i^\epsilon(1+\frac{\bar{u}_i^\epsilon}{m_i^\epsilon})^2
       +2\pi d_i^\epsilon(1-\frac{\bar{u}_i^\epsilon}{m_i^\epsilon})+
    \pi(1-\frac{\bar{u}_2^\epsilon}{m_1^\epsilon})^2(\log{L^2}-\log{\delta^2})\right)\\[1.7ex]
        \s\s +\sum\limits_{i=1,2}\frac{A'\bar{u}_i^\epsilon}
           {(m_i^\epsilon)^2}+o(1)+O(\delta^2).
\end{array}$$
We set $s_i^\epsilon=1+\frac{\bar{u}_i^\epsilon}{m_i^\epsilon}$. Then, for fixed $L$, $\delta$, we have
$$\frac{1}{3}\dint_{\Sigma}(|\nabla u_1|^2+|\nabla u_2|^2
+\nabla u_1\nabla
u_2)dV_g+(4\pi-\epsilon)\bar{u}_1^\epsilon+(4\pi-\epsilon)\bar{u}_2^\epsilon\geq
\sum_im_i^\epsilon(s_i^\epsilon+O(\frac{1}{m_i^\epsilon}))^2 +C.$$
Since $\Phi_\epsilon(u_\epsilon)\leq C$, we see that
$$|s_i^\epsilon|=O(\frac{1}{m_i^\epsilon}).$$
Hence for both $i=1,2$, $s_i^\epsilon\rightarrow 0$ as
$\epsilon\to 0$. So,
$$\begin{array}{l}
\frac{1}{3}\dint_{B_\delta(x_1^\epsilon)\cup B_\delta(x_2^\epsilon)
}(|\nabla u_1|^2+|\nabla u_2|^2
+\nabla u_1\nabla u_2)dV_g+(4\pi-\epsilon)\bar{u}_1^\epsilon+(4\pi-\epsilon)
\bar{u}_2^\epsilon\\[1.7ex]
\s\geq \frac{1}{2}\dint_{B_L}|\nabla w|^2dx+4\pi d_1^\epsilon+4\pi d_2^\epsilon
+8\pi(\log{L^2}-\log{\delta^2})+o(1)+O(\delta^2)\\[1.7ex]
\s=
\frac{1}{2}\dint_{B_L}|\nabla w|^2dx+8\pi w(L)+8\pi(\log{L^2}-\log{\delta^2})\\[1.7ex]
\s\s-4\pi\sup_{\partial B_\delta(p_1)}G_1
-4\pi\sup_{\partial B_\delta(p_2)}G_2+o(1)+O(\delta^2).
\end{array}\eqno(3.1)$$
By a direct calculation, we obtain
$$\dint_{B_L}|\nabla w|^2dx=16\pi\log{(1+\pi L^2)}-\frac{16\pi^2 L^2}{1+\pi L^2}.\eqno(3.2)$$
Moreover, by (2.1) and (2.2),  we have
$$\begin{array}{l}\frac{1}{3}\dint_{B^c_\delta(x_1^\epsilon)\cap B^c_\delta(x_2^\epsilon)
}(|\nabla u_1|^2+|\nabla u_2|^2
+\nabla u_1\nabla u_2)dV_g\\[1.7ex]
\s =\frac{1}{3}\dint_{B^c_\delta(x_1^\epsilon)\cap B^c_\delta(x_2^\epsilon)
}(|\nabla G_1|^2+|\nabla G_2|^2
+\nabla G_1\nabla G_2)dV_g+o(1)\\[1.7ex]
\s =-\frac{1}{3}\sum\limits_{i=1,2}\dint_{\partial B_\delta(p_i)}(G_1
\frac{\partial G_1}{\partial n}+G_2
\frac{\partial G_2}{\partial n}+\frac{G_1
\frac{\partial G_2}{\partial n}+G_2
\frac{\partial G_1}{\partial n}}{2})dS_g+o(1)\\[1.7ex]
\s\s\s +\dint_{B_{\delta}(p_1)+B_{\delta}(p_2)}2\pi(G_1+G_2)dV_g+o(1)\\[1.7ex]
\s =-\frac{1}{3}\sum\limits_{i=1,2}\dint_0^{2\pi}(G_1
\frac{\partial G_1}{\partial r}+G_2
\frac{\partial G_2}{\partial r}+\frac{G_1
\frac{\partial G_2}{\partial r}+G_2
\frac{\partial G_1}{\partial r}}{2})rd\theta|_{r=\delta}+o(1)\\[1.7ex]
\s\s\s +\dint_{B_{\delta}(p_1)+B_{\delta}(p_2)}2\pi(G_1+G_2)dV_g+o(1)\\[1.7ex]
\s =-16\pi\log{\delta}-2\pi A_1(p_1)-2\pi A_2(p_2)+o(1)+O(\delta\log \delta).
\end{array}\eqno(3.3)$$
In the end, (3.1), (3.2) and (3.3) imply that
$$\inf\Phi_0(u)\geq -8\pi\log{\pi}-8\pi-2\pi(A_1(p_1)+A_2(p_2)).$$

\section{Lower bound for case 2}
In this case, we
set $v_2^\epsilon=\frac{1}{3}(2u_2^\epsilon+u_1^\epsilon)-\frac{1}{3}
(2\bar{u}_2^\epsilon+\bar{u}_1^\epsilon)$. Then, we have
$$\left\{\begin{array}{l}
             -\Delta v_2^\epsilon=(4\pi-\epsilon)e^{u_2^\epsilon}-(4\pi-\epsilon)\\[1.7ex]
             \dint v_2^\epsilon=0
         \end{array}\right.$$
By the standard elliptic estimates,  $||v_2||_{C^1(M)}<C$.

Similar to the case 1, we have
$$\frac{1}{3}\dint_{B_\delta(x_1^\epsilon)}(|\nabla u_1|^2+|\nabla u_2|^2
+\nabla u_1\nabla u_2)dV_g=
\frac{1}{4}\dint_{B_\delta(x_1^\epsilon)}|\nabla u_1^\epsilon|^2dV_g+O(\delta^2),$$
and
$$\begin{array}{l}
\frac{1}{4}\dint_{B_{\delta}(x_1^\epsilon)}|\nabla u_1|^2dV_g+(4\pi-\epsilon)
\bar{u}_1^\epsilon\\[1.7ex]
\s \geq
\frac{1}{4}\dint_{B_L}|\nabla w|^2dx+
\left(\pi\frac{(m_1^\epsilon+\bar{u}_1^\epsilon)^2}{m_1^\epsilon}
+2\pi d_1^\epsilon(1-\frac{\bar{u}_1^\epsilon}{m_1^\epsilon})+
\pi(1-\frac{\bar{u}_1^\epsilon}{m_1^\epsilon})^2(\log{L^2}-\log{\delta^2})\right)\\[1.7ex]
\s\s +\frac{A'\bar{u}_1^\epsilon}
           {(m_1^\epsilon)^2}+o(1)+O(\delta^2).
\end{array}$$
By an argument similar to the one used in the case 1, we can show
that $\frac{\bar{u}_1^\epsilon}{m_1^\epsilon}\rightarrow-1$, hence
$$\begin{array}{l}
\frac{1}{3}\dint_{B_\delta(x_1^\epsilon)}(|\nabla u_1|^2+|\nabla u_2|^2
+\nabla u_1\nabla u_2)dV_g+(4\pi-\epsilon) \bar{u}_1^\epsilon\\[1.7ex]
\s\geq
\frac{1}{4}\dint_{B_L}|\nabla w|^2dx+4\pi w(L)+4\pi(\log{L^2}-\log{\delta^2})
-4\pi\sup_{\partial B_\delta(p_1)}G_1
+o(1)+O(\delta^2).
\end{array}$$
Set
$$G_1=-4\log{r}+A_1(p)+o(x),\s G_2=2\log{r}+A_2(p)+o(x).$$
Applying (2.4), we get
$$\begin{array}{l}\frac{1}{3}\dint_{B^c_\delta(x_1^\epsilon)
}(|\nabla u_1|^2+|\nabla u_2|^2
+\nabla u_1\nabla u_2)dV_g\\[1.7ex]
\s =\frac{1}{3}\dint_{B^c_\delta(x_1^\epsilon)
}(|\nabla G_1|^2+|\nabla G_2|^2
+\frac{\nabla G_1\nabla G_2+\nabla G_2\nabla G_1}{2})dV_g+o(1)\\[1.7ex]
\s =\dint_{\partial B_\delta(p)}(G_1
\frac{\partial G_1}{\partial n}+G_2
\frac{\partial G_2}{\partial n}+\frac{G_1
\frac{\partial G_2}{\partial n}+G_2
\frac{\partial G_1}{\partial n}}{2})\\[1.7ex]
\s +2\pi\dint_{B_\delta(p)}(G_1+G_2)dV_g-
2\pi\dint_\Sigma G_2dV_g+o(1)+O(\delta\log \delta)\\[1.7ex]
\s =-8\pi\log{\delta}-2\pi A_1(p_1)-2\pi\dint_\Sigma G_2dV_g+o(1)+O(\delta\log\delta).
\end{array}$$
In the end, we obtain
$$\inf\Phi_0(u)\geq -4\pi\log{\pi}-
2\pi A_1(p)+2\pi\dint G_2dV_g.$$

\section{Test functions for case 1}
In this section, we will construct a function
$\phi=(\phi_1,\phi_2)\in H^{1,2}(M)\times H^{1,2}(M)$, such that
$$\Phi_0(\phi)<-8\pi\log{\pi}-8\pi-2\pi(A_2(p_1)+A_1(p_2)),$$
whenever (1.4) holds. So, under the assumption (1.4), the case 1
will not happen.

Let $(\Omega_i;(x,y))$ be an isothermal coordinate system around
$p_i$ ($i=1,2$). We set
$$r(x,y)=\sqrt{x^2+y^2},\s and\s B_\delta=\{(x,y):x^2+y^2<\delta^2\}.$$
We assume that near $p_i$ ($i=1,2$), for each $k=1,2$,
\begin{eqnarray*}G_k&=&a_k(p_i)\log{r}+A_k(p_i)+\lambda_k(p_i)x+\mu_k(p_i)y
\\&&+\alpha_k(p_i)x^2+\beta_k(p_i)y^2
+\gamma_k(p_i)xy+h(x,y)+O(r^4).\end{eqnarray*} We have
$a_1(p_1)=a_2(p_2)=-4$, and $a_1(p_2)=a_2(p_1)=2$. Moreover, we
assume that
$$g|_{\Omega_i}=e^{\varphi_i}(dx^2+dy^2),$$
and
$$\varphi_i=b_1(p_i)x+b_2(p_i)y
+c_1(p_i)x^2+c_2(p_i)y^2+c_{12}(p_i)xy+O(r^3).$$
It is well known that
$$K(p_i)=-(c_1(p_i)+c_2(p_i)),$$
$$|\nabla u|^2dV_g=|\nabla u|^2dxdy,$$
and
$$\frac{\partial u}{\partial n}dS_g=\frac{\partial u}{\partial r}rd\theta,\s(S=\partial B_r).$$
For $\alpha_k$ and $\beta_k$, we have the following lemma:
\bthm{Lemma}{5.1} For any $k,i$, we have
$$\alpha_k(p_i)+\beta_k(p_i)=2\pi.$$\ethm
\proof Near $p_i$, we have
$$2\alpha_k(p_i)+2\beta_k(p_i)+O(r)=\Delta_0 G_k(x,y)=e^{-\varphi_i}4\pi.$$
\begin{flushright}
{\small $\Box$}
\end{flushright}
\hspace{2ex}

We choose
$$\phi_1=\left\{
     \begin{array}{ll}
       w(\frac{x}{\epsilon})+\lambda_1(p_1)r\cos\theta+
        \mu_1(p_1)r\sin\theta&(x,y)\in B_{L\epsilon}(p_1)\\[1.7ex]
       G_1-\eta_1 H^{p_1}_1+4\log{L\epsilon}-2\log{(1+\pi L^2)}
        -A_1(p_1)&(x,y)\in B_{2L\epsilon}\setminus B_{L\epsilon}(p_1)\\[1.7ex]
       G_1-\eta_2 H^{p_2}_1+4\log{L\epsilon}-2\log{(1+\pi L^2)}-A_1(p_1)&(x,y)\in
        B_{2L\epsilon}\setminus B_{L\epsilon}(p_2)\\[1.7ex]
       -\frac{\omega(\frac{x}{\epsilon})+2\log{(1+\pi L^2)}}{2}
        +\lambda_1(p_2)r\cos\theta+\mu_1(p_2)r\sin\theta&\\
       \s\s\s +6\log{L\epsilon}-2\log{(1+\pi L^2)}+A_1(p_2)-A_1(p_1)&(x,y)\in
        B_{L\epsilon}(p_2)\\[1.7ex]
              G_1+4\log{L\epsilon}-2\log{(1+\pi L^2)}-A_1(p_1)&others,
      \end{array}\right.$$
and
$$\phi_2=\left\{
     \begin{array}{ll}
       w(\frac{x}{\epsilon})+\lambda_2(p_2)r\cos\theta+\mu_2(p_2)r\sin\theta &(x,y)\in
       B_{L\epsilon}(p_2)\\[1.7ex]
       G_2-\eta_2 H^{p_2}_2+4\log{L\epsilon}-2\log{(1+\pi L^2)}
        -A_2(p_2)&(x,y)\in B_{2L\epsilon}\setminus B_{L\epsilon}(p_2)\\[1.7ex]
       G_2-\eta_1 H^{p_1}_2+4\log{L\epsilon}-2\log{(1+\pi L^2)}-A_2(p_2)&(x,y)\in
        B_{2L\epsilon}\setminus B_{L\epsilon}(p_1)\\[1.7ex]
       -\frac{\omega(\frac{x}{\epsilon})+2\log{(1+\pi L^2)}}{2}+
       \lambda_2(p_1)r\cos\theta+\mu_2(p_1)r\sin\theta\\
       \s\s\s +
       6\log{L\epsilon}-2\log{(1+\pi L^2)}+A_2(p_1)-A_2(p_2)&(x,y)\in
        B_{L\epsilon}(p_1)\\[1.7ex]
              G_2+4\log{L\epsilon}-2\log{(1+\pi L^2)}-A_2(p_2)&others.
      \end{array}\right.$$
Here,
$$H^{p_i}_k=G_k-a_k(p_i)\log r-A_k(p_i)-\lambda_k(p_i)r\cos\theta-\mu_k(p_i)r\sin\theta,$$
and $\eta_i$ is a cut-off function which equals 1 in
$B_{L\epsilon}(p_i)$, equals $0$ in $B_{2L\epsilon}^c(p_i)$. We
may assume that
$$|\nabla\eta_i|\leq \frac{1}{L\epsilon}.$$

Now, we compute $\Phi_0(\phi)$.

Firstly, we compute $\int_\Sigma|\nabla\phi_1|^2dV_g$ and $\int_\Sigma|\nabla\phi_2|^2dV_g$.

Let $\Omega=\Sigma\setminus (B_{L\epsilon}(p_1)\cup
B_{L\epsilon}(p_2))$. Then
$$\begin{array}{ll}
   \dint_\Sigma|\nabla\phi_1|^2dV_g=&\dint_{B_{L\epsilon}(p_1)\cup
       B_{L\epsilon}(p_2)}|\nabla\phi_1|^2dxdy
       +\dint_\Omega|\nabla G_1|^2dV_g-\\[1.7ex]
     & 2\sum\limits_{i=1,2}\dint_{\Sigma}\nabla G_1\nabla\eta_iH^{p_i}_1dV_g
       +\sum\limits_{i=1,2}\dint_{\Sigma}|\nabla\eta_iH^{p_i}_1|^2dV_g.
\end{array}$$
It is clear that we have
$$\dint_{B_{L\epsilon}(p_2)}|\nabla\phi_1|^2dV_g=\frac{1}{4}\dint_{B_L}|\nabla w|^2dxdy+
\pi(L\epsilon)^2
(\lambda_1^2(p_2)+\mu_1^2(p_2)),$$
and
$$\dint_{B_{L\epsilon}(p_1)}|\nabla\phi_1|^2dV_g=
\dint_{B_L}|\nabla w|^2dxdy+\pi(L\epsilon)^2(
\lambda_1^2(p_1)+\mu_1^2(p_1)).$$
Calculating directly and using
the fact that $\int_0^{2\pi}hd\theta=0$, we obtain,
$$\begin{array}{ll}
    \dint_\Sigma \nabla G_1\nabla\eta_1 H^{p_1}_1dV_g&=-\dint_{\partial B_{L\epsilon}(p_1)}
                  \frac{\partial G_1}{\partial n}H_1^{p_1}dS_g-4\pi\dint_{B_{2L\epsilon
                   \setminus B_{L\epsilon}}}\eta_1H_1^{p_1}dV_g\\[1.7ex]
               &=-\dint_0^{2\pi}(-\frac{4}{r}+\lambda_1(p_1)\cos\theta+
                 \mu_1(p_1)\sin\theta+O(r))\\[1.7ex]
                 &\s\s\s \times(\alpha_1(p_1)r^2
                 \cos^2\theta+\beta_1(p_1)r^2\sin^2\theta+h+O(r^4))rd\theta\\[1.7ex]
               &=4\pi(\alpha_1(p_1)+\beta_1(p_1))(L\epsilon)^2+O(L\epsilon)^4\\[1.7ex]
               &=8\pi^2(L\epsilon)^2+O(L\epsilon)^4.
\end{array}$$
Similarly, we get
$$\dint_\Sigma \nabla G_1\nabla\eta_2 H^{p_2}_1dV_g=-4\pi^2(L\epsilon)^2+O(L\epsilon)^4.$$
It is obvious that
$$\dint_\Sigma|\nabla\eta_jH^{p_j}_i|^2dV_g=\dint_{B_{2L\epsilon}\setminus B_{L\epsilon}(p_j)}O(r^2)dV_g
=O((L\epsilon)^4).$$ Hence
\begin{eqnarray*}\dint_\Sigma|\nabla\phi_1|^2dV_g&=&\frac{5}{4}\dint_{B_L}|\nabla
w|^2dx- 8\pi^2(L\epsilon)^2+\dint_\Omega|\nabla G_1|^2dV_g\\&&+
\pi(L\epsilon)^2\sum_{i=1,2}(\lambda^2_1(p_i)+\mu_1^2(p_i))+O(L\epsilon)^4.\end{eqnarray*}

In the same way, we can show that
\begin{eqnarray*}\dint_\Sigma|\nabla\phi_2|^2dV_g&=&\frac{5}{4}\dint_{B_L}|\nabla w|^2dx-
8\pi^2(L\epsilon)^2+\dint_\Omega|\nabla G_2|^2dV_g\\&&+
\pi(L\epsilon)^2\sum_{i=1,2}(\lambda^2_2(p_i)+\mu_2^2(p_i))+O(L\epsilon)^4.\end{eqnarray*}

Next, we compute $\int_\Sigma\nabla\phi_1\nabla\phi_2dV_g$. We have
$$\begin{array}{ll}
   \dint_\Sigma\nabla\phi_1\nabla\phi_2dV_g&=
     \sum\limits_{i=1,2}\dint_{B_{L\epsilon}(p_i)}
       \nabla\phi_1\nabla\phi_2dV_g+\dint_\Omega\nabla G_1\nabla G_2dV_g\\[1.7ex]
    &\s-\dint_\Sigma\nabla G_1\nabla\eta_1 H^{p_1}_2dV_g-
      \dint_\Sigma \nabla G_2\nabla\eta_1 H^{p_1}_1dV_g-
      \dint_\Sigma\nabla G_1\nabla\eta_2 H^{p_2}_2dV_g\\[1.7ex]
    &-\dint_\Sigma\nabla G_2\nabla\eta_2H^{p_2}_1dV_g\s+
     \s\sum\limits_{i=1,2}\dint_\Sigma\nabla\eta_iH^{p_i}_1\nabla\eta_iH^{p_i}_2dV_g\\[1.7ex]
    &=-\dint_{B_L}|\nabla w|^2dx+\pi(L\epsilon)^2\sum_{i=1,2}(\lambda_1(p_i)\lambda_2(p_i)
       +\mu_1(p_i)\mu_2(p_i))\\[1.7ex]
    &\s\s\s-8\pi(L\epsilon)^2+\dint_\Omega\nabla G_1\nabla G_2dV_g+O((L\epsilon)^4).
\end{array}$$

Then, we calculate $\dint_\Omega(|\nabla G_1|^2+|\nabla
G_2|^2+\nabla G_1\nabla G_2) dV_g$. We have
\begin{eqnarray*}
    \lefteqn{\dint_\Omega(|\nabla G_1|^2+|\nabla G_2|^2+\nabla G_1\nabla
    G_2)dV_g}\\
        &=&\int_\Omega(|\nabla G_1|^2+|\nabla G_2|^2+\frac{\nabla G_1\nabla G_2+
          \nabla G_2\nabla G_1}{2})dV_g\\&=&-\dint_{\partial B_{L\epsilon}(p_1)+\partial B_{L\epsilon}(p_2)}
        (G_1\frac{\partial G_1}{\partial n}+G_2\frac{\partial G_2}{\partial n}
         +\frac{G_1\frac{\partial G_2}{\partial n}+
          G_2\frac{\partial G_1}{\partial n}}{2})dS_g\\
   && +6\pi\dint_{B_{L\epsilon(p_1)}+B_{L\epsilon}(p_2)}(G_1+G_2)dV_g
\end{eqnarray*}

\bthm{Lemma}{5.2} For any $k,m,i=1,2$, we have
\begin{eqnarray*}
\dint_{\partial B_{r}(p_i)}G_k\frac{\partial G_m}{\partial
n}dS_g&=& 2\pi a_k(p_i)a_m(p_i)\log{r}+2\pi^2a_k(p_i)r^2\\&&+\pi
r^2(
\lambda_k(p_i)\lambda_m(p_i)+\mu_k(p_i)\mu_m(p_i))\\
&& +2\pi a_k(p_i)A_m(p_i)+
4\pi^2r^2A_m(p_i)\\&&+4\pi^2r^2a_m(p_i)\log{r}+O(r^4\log r).
\end{eqnarray*}
\proof Since $\int_0^{2\pi}h(r,\theta)d\theta=\int_0^{2\pi}\frac{\partial h}{\partial r}(r,\theta)d\theta=0$, we have
\begin{eqnarray*}
\dint_{\partial B_{r}(p_i)}G_k\frac{\partial G_m}{\partial
n}dS_g&=&
  \dint_0^{2\pi}(\frac{a_k(p_i)}{r}+\lambda_k(p_i)\cos\theta+\mu_k(p_i)\sin\theta\\&&+
   2r\alpha_k(p_i)\cos^2\theta+2r\beta_k(p_i)\sin^2\theta+2r\gamma_k(p_i)\sin\theta
   \cos\theta)\\&& \times(a_m(p_i)\log{r}+A_m(p_i)+\lambda_m(p_i)r\cos\theta
   +\mu_m(p_i)r\sin\theta\\&& +r^2\alpha_m(p_i)\cos^2\theta+r^2\beta_m(p_i)
   \sin^2\theta+r^2\gamma_m(p_i)\sin\theta
   \cos\theta) rd\theta\\&&+O(r^4\log{r})\\
  &=&2\pi a_k(p_i)a_m(p_i)\log{r}+\pi(\alpha_m(p_i)+\beta_m(p_i))a_k(p_i)r^2\\&&
   +\pi r^2(\lambda_k(p_i)\lambda_m(p_i)+\mu_k(p_i)\mu_m(p_i))
     +2\pi a_k(p_i)A_m(p_i)\\&& +
     (2\pi r^2A_m(p_i)+2\pi r^2a_m(p_i)\log{r})(\alpha_k(p_i)+\beta_k(p_i))\\&&+O(r^4\log r).
\end{eqnarray*}
\begin{flushright}
{\small $\Box$}
\end{flushright}
\hspace{2ex}

Then
$$\begin{array}{ll}
\dint_{\partial B_r(p_1)}G_1\frac{\partial G_1}{\partial n}dS_g=&
32\pi \log{r}-8\pi^2r^2+\pi r^2(
\lambda^2_1(p_1)+\mu_1^2(p_1))\\[1.7ex]
& -8\pi A_1(p_1)+
4\pi^2r^2A_1(p_1)-16\pi^2r^2\log{r}+O(r^4\log r).
\end{array}$$

$$\begin{array}{ll}
\dint_{\partial B_r(p_1)}G_2\frac{\partial G_2}{\partial n}dS_g=&
8\pi\log{r}+4\pi^2r^2+\pi r^2(
\lambda_2^2(p_1)+\mu_2^2(p_1))\\[1.7ex]
& +4\pi A_2(p_1)+
4\pi^2r^2A_2(p_1)+8\pi^2r^2\log{r}+O(r^4\log r).
\end{array}$$

$$\begin{array}{ll}
\dint_{\partial B_r(p_1)}G_1\frac{\partial G_2}{\partial n}dS_g=&
-16\pi \log{r}-8\pi^2r^2+\pi r^2(
\lambda_1(p_1)\lambda_2(p_1)+\mu_1(p_1)\mu_2(p_1))\\[1.7ex]
& -8\pi A_2(p_1)+
4\pi^2r^2A_2(p_1)+8\pi^2r^2\log{r}+O(r^4\log r).
\end{array}$$

$$\begin{array}{ll}
\dint_{\partial B_r(p_1)}G_2\frac{\partial G_1}{\partial n}dS_g=&
-16\pi \log{r}+4\pi^2r^2+\pi r^2(
\lambda_2(p_1)\lambda_1(p_1)+\mu_2(p_1)\mu_1(p_1))\\[1.7ex]
& +4\pi A_1(p_1)+
4\pi^2r^2A_1(p_1)-16\pi^2r^2\log{r}+O(r^4\log r).
\end{array}$$

Hence
$$\begin{array}{l}
\dint_\Omega(|\nabla G_1|^2+|\nabla G_2|^2)dV_g
=-(80\pi\log {L\epsilon}-8\pi^2(L\epsilon)^2+\pi (L\epsilon)^2\sum_{i,j=1,2}(\lambda_i^2(j)+\mu_i^2(j))\\[1.7ex]
\s\s\s -8\pi A_1(p_1)-8\pi A_2(p_2)+4\pi^2(L\epsilon)^2(A_1(p_1)+A_2(p_2)+A_2(p_1)+A_1(p_2))\\[1.7ex]
\s\s\s-
16\pi^2(L\epsilon)^2\log {L\epsilon}+4\pi A_2(p_1)+4\pi A_1(p_2))+4\pi\dint_{\Omega}G_1dV_g+
4\pi\dint_{\Omega}G_2dV_g\\[1.7ex]
\s\s\s+O((L\epsilon)^4\log {L\epsilon}),
\end{array}$$
and
$$\begin{array}{l}
\dint_{\Omega}\nabla G_1\nabla G_2dV_g=-(-32\log{{L\epsilon}}-4\pi^2(L\epsilon)^2+\pi (L\epsilon)^2\sum_{i\neq j}
(\lambda_i(p_j)\lambda_j(p_i)+\mu_i(p_j)\mu_j(p_i))\\[1.7ex]
\s\s\s-4\pi A_2(p_1)-4\pi A_1(p_2)+2\pi A_1(p_1)+2\pi A_2(p_2)\\[1.7ex]
\s\s\s+2\pi^2(L\epsilon)^2(A_1(p_1)+A_1(p_2)+A_2(p_1)+A_2(p_2))-8\pi^2(L\epsilon)^2\log {L\epsilon})\\[1.7ex]
\s\s\s+2\pi\sum_{i=1,2}\dint_{B_{L\epsilon}(p_i)}(G_1+G_2)dV_g+O((L\epsilon)^4\log {L\epsilon}).
\end{array}$$
It is easy to check that
$$\begin{array}{lll}
 6\pi\dint_{B_{L\epsilon}(p_1)+B_{L\epsilon}(p_2)
  }(G_1+G_2)dV_g&=&6\pi^2 (L\epsilon)^2(\sum_{i,j=1,2} A_i(p_j))
   -24\pi^2(L\epsilon)^2\log {L\epsilon}\\
   &&+12\pi^2 (L\epsilon)^2+O((L\epsilon)^4\log {L\epsilon}).
  \end{array}$$
So, we get
\begin{eqnarray}\dint_{\Sigma}(|\nabla\phi_1|^2+|\nabla\phi_2|^2+\nabla\phi_1\nabla\phi_2)dV_g
&=&\frac{3}{2}\dint_{B_L}|\nabla w|^2dxdy-48\pi\log {L\epsilon}+6\pi
A_1(p_1)\nonumber\\&&+6\pi A_2(p_2)+O((L\epsilon)^4\log{{L\epsilon}}).
\end{eqnarray}

We calculate $\dint_\Sigma(\phi_1+\phi_2)dV_g$. We have
$$\begin{array}{ll}
    \dint_\Sigma\phi_1dV_g&=\epsilon^2\dint_{B_L}we^{\varphi(\epsilon x,\epsilon y)}dxdy
    -\dint_\Sigma(\eta_1 H^{p_1}_1+\eta_2 H^{p_2}_1)dV_g\\[1.7ex]
    &\s+(4\log{{L\epsilon}}-2\log{(1+\pi L^2)}-A_1(p_1))(1-\dint_{B_{L\epsilon}(p_1)}dV_g)+2\log {L\epsilon}\dint_{B_{L\epsilon}(p_2)}
    dV_g\\[1.7ex]
    &\s+A_1(p_2)\dint_{B_{L\epsilon}(p_2)}dV_g-\dint_{B_{L\epsilon}(p_1)+B_{L\epsilon}(p_2)}G_1
    -\epsilon^2\dint_{B_L}\frac{w+2\log{(1+\pi L^2)}}{2}dV_g\\[1.7ex]
    &\s+\sum\limits_{i=1,2}\dint_{B_{L\epsilon}(p_i)}
      (\lambda_1(p_i)x+\mu_1(p_i)y)dV_g.
\end{array}$$
Since
$$\begin{array}{ll}
\dint_{B_{L\epsilon}(p_1)+B_{L\epsilon}(p_2)}G_1dV_g&=\dint_0^{L\epsilon}
 (-2\log{r}+A_1(p_1)+A_1(p_2))2\pi rdr+O((L\epsilon)^4\log{L\epsilon})\\[1.7ex]
 &=-2\pi (L\epsilon)^2\log{{L\epsilon}}+
   \pi (L\epsilon)^2+(A_1(p_1)+A_1(p_2))\pi (L\epsilon)^2+O((L\epsilon)^4\log{L\epsilon}),
\end{array}$$
we have
\begin{eqnarray*}\dint_\Sigma\phi_1dV_g&=&\frac{\epsilon^2}{2}\dint_{B_L}we^{\varphi(\epsilon
x,\epsilon y)}dxdy
 +4\log{{L\epsilon}}+\pi (L\epsilon)^2\log{(1+\pi L^2)}\\&&
    -\pi (L\epsilon)^2-A_1(p_1)-2\log(1+\pi L^2)
    +O((L\epsilon)^4\log {L\epsilon}).\end{eqnarray*}
Similarly, we have
\begin{eqnarray*}\dint_\Sigma\phi_2dV_g&=&\frac{\epsilon^2}{2}\dint_{B_L}we^{\varphi(\epsilon x,\epsilon y)}dxdy
+4\log{{L\epsilon}}+\pi (L\epsilon)^2\log{(1+\pi L^2)})\\&&
    -\pi (L\epsilon)^2-A_2(p_2)-2\log(1+\pi L^2)
    +O((L\epsilon)^4\log {L\epsilon}).\end{eqnarray*}
Moreover, we have
$$\dint_{B_L}\omega e^{\varphi(\epsilon x,\epsilon y)}dxdy
=2\pi L^2-2\log(1+\pi L^2)-2\pi L^2\log(1+\pi L^2)+O(L^2\epsilon^2\log{L}),$$
hence, we get
\begin{eqnarray}\dint_\Sigma(\phi_1+\phi_2)dV_g
    &=&-A_1(p_1)-A_2(p_2)+8\log{{L\epsilon}}-4\log(1+\pi L^2)\nonumber\\&&-2\epsilon^2\log(1+\pi L^2)
    +O((L\epsilon)^4\log {L\epsilon}).\end{eqnarray}

We denote $B(p_j)=\frac{(b_1(p_j)+\lambda_1(p_j))^2
+(b_2(p_j)+\lambda_2(p_j))^2}{4}$, and $M_i=\frac{-\frac{K(p_i)}{2}+B(p_i)}{\pi}$.
Then, we have
\begin{eqnarray*}
   \lefteqn{\dint_{B_{L\epsilon}(p_1)}e^{\phi_1}dV_g}\\&=&\epsilon^2\dint_{B_L}
        \frac{e^{(b_1(p_1)+\lambda_1(p_1))\epsilon x+(b_2(p_1)+\lambda_2(p_1))
         \epsilon y+
          c_1(p_1)\epsilon^2x^2+c_2(p_1)\epsilon^2y^2+c_{12}\epsilon^2xy+O((r\epsilon)^3)}
             }{(1+\pi r^2)^2}dxdy\\
       &=&\epsilon^2\dint_0^L\frac{1+\epsilon^2\pi M_ir^2+O(\epsilon^3 r^3)
             }{(1+\pi r^2)^2}2\pi r dr\\
       &=&\epsilon^2(1-\epsilon^2 M_i)
         \frac{\pi L^2}{1+\pi L^2}
         +\epsilon^4 M_i\log(1+\pi L^2)+\epsilon^2O(\epsilon^3\log{L})\\
       &=&\epsilon^2\left(1-\frac{1}{1+\pi L^2}+\epsilon^2M_i\log(1+\pi L^2)+
        O(\epsilon^2)+O(\epsilon^3\log L)\right),
\end{eqnarray*}
and we also have
\begin{eqnarray*}
 \lefteqn{\dint_{B_\delta(p_1)\setminus
 B_{L\epsilon(p_1)}}e^{\phi_1}dV_g}\\&=&
    \frac{(L\epsilon)^4}{(1+\pi L^2)^2}\dint_{L\epsilon}^{\delta}e^{-4\log r+(\lambda_1(p_1)+b_1(p_1))x+(
    \mu_1(p_1)+b_2(p_1))y}\\
  &&\times e^{(c_1(p_1)+\alpha_1(p_1))x^2+(c_2(p_1)+\beta_1(p_1))y^2+(c_{12}(p_1)+\gamma_1)xy+(1-\eta_1)H_1^{p_1}+O(r^3)}2\pi rdr\\
  &=&\epsilon^2\left(\frac{\pi L^2}{(1+\pi L^2)^2}-2(M_1+1)\epsilon^2\log{L\epsilon}+O(\epsilon^2)+
  O(\frac{1}{L^4})\right).
\end{eqnarray*}

Since outsider $B_\delta(p_1)$, $G_1$ is bounded above, we have
$$\dint_{\Sigma\setminus B_\delta}e^{\phi_1}=O(\epsilon^4).$$
Notice that $\frac{\pi L^2}{(1+\pi L^2)^2}-\frac{1}{1+\pi
L^2}=O(\frac{1}{L^4})$, we get
\begin{eqnarray}\log\dint_\Sigma
e^{\phi_1}dV_g&=&\log{\epsilon^2} +\epsilon^2M_1\log(1+\pi L^2)-2\epsilon^2(M_1+1)\log{L\epsilon}\nonumber\\
&&O(\epsilon^2)+O(\frac{1}{L^4}).\end{eqnarray}

In the same way, we can get
\begin{eqnarray}\log\dint_\Sigma e^{\phi_2}dV_g&=&\log{\epsilon^2}+
\epsilon^2M_2\log(1+\pi L^2)-2\epsilon^2(M_2+1)\log{L\epsilon}\nonumber\\
&&+ O(\epsilon^2)+O(\frac{1}{L^4}).\end{eqnarray}

It follows from (5.1), (5.2), (5.3) and (5.4) that
\begin{eqnarray*}
    \Phi_0(\phi)&=&\frac{1}{2}\dint_{B_L}|\nabla w|^2dxdy+16\log{L\epsilon} -2\pi(A_1(p_1)+A_2(p_2))
     \\&&-16\pi\log(1+\pi L^2)
     -8\pi\log\epsilon^2 -8\pi\epsilon^2\log(1+\pi L^2)\\
     &&-4\pi\epsilon^2\left((M_1+M_2)\log(1+\pi L^2)-2(M_1+M_2+2)\log{L\epsilon}\right)
     \\&&+O(\frac{1}{L^4})+O(\epsilon^2)+O((L\epsilon)^4\log{L\epsilon})+O(\epsilon^3
    \log{L})\\
    &=&-8\pi\log\frac{1+\pi L^2}{ L^2}-8\pi\frac{\pi L^2}{1+\pi L^2}
     -2\pi(A_1(p_1)+A_2(p_2))\\&&
     -4\pi\epsilon^2(M_1+M_2+2)(\log(1+\pi L^2)-2\log{L\epsilon})\\&&
     +O(\frac{1}{L^4})+O(\epsilon^2)+O((L\epsilon)^4\log{L\epsilon})+O(\epsilon^3
    \log{L})\\
    &=&-8\pi\log\pi-8\pi-2\pi(A_1(p_1)+A_2(p_2))\\
    &&-4\pi(M_1+M_2+2)\epsilon^2
    (\log(1+\pi L^2)-2\log{L\epsilon})\\&&+O(\frac{1}{L^4})+O(\epsilon^2)+O((L\epsilon)^4\log{L\epsilon})+O(\epsilon^3
    \log{L}).
\end{eqnarray*}
Under the assumption (1.1), we have $M_1+M_2+2>0$. Let
$L^4\epsilon^2=\frac{1}{\log(-\log\epsilon)}$, we get
\begin{eqnarray*}
    \Phi_0(\phi)
    &=&-8\pi\log\pi-8\pi-2\pi(A_1(p_1)+A_2(p_2))\\
    &&-4\pi(M_1+M_2+2)\epsilon^2
    (-\log{\epsilon^2})+o(\epsilon^2(-\log\epsilon^2)).
\end{eqnarray*}
Then for
sufficiently small $\epsilon$, we have
$$\Phi_0(\phi)<-8\pi\log\pi-8\pi-2\pi(A_1(p_1)+A_2(p_2)).$$ This
proves our claim.

\section{Test functions for case 2}
Assume that (1.4) holds on $\Sigma$, we will construct a function
$\phi=(\phi_1,\phi_2)\in H^{1,2}(M)\times H^{1,2}(M)$, such that
$$\Phi_0(\phi)<-4\pi\log{\pi}-
2\pi A_1(p)+2\pi\dint G_2dV_g.$$
Let $(\Omega;(x,y))$ be an isothermal  coordinate system around $p$.
We assume that near $p$
$$G_k=a_k\log{r}+A_k(p)+\lambda_kx+\mu_ky
+\alpha_kx^2+\beta_ky^2+\gamma_kxy+h_k(x,y)+O(r^4).$$
We have $a_1(p)=-4$, and $a_2(p_1)=2$. Moreover, we assume that
$$g|_{\Omega}=e^{\varphi}(dx^2+dy^2),$$
and
$$\varphi=b_1x+b_2y
+c_1x^2+c_2y^2+c_{12}xy+O(r^3).$$ Similar to the case 1, there
hold
$$\alpha_k+\beta_k=2\pi,\s i=1,2.$$

We choose
$$\phi_1=\left\{
     \begin{array}{ll}
       w(\frac{x}{\epsilon})+\lambda_1r\cos\theta+
        \mu_1r\sin\theta&x\in B_{L\epsilon}(p)\\[1.7ex]
       G_1-\eta H_1+4\log{L\epsilon}-2\log{(1+\pi L^2)}
        -A_1(p)&x\in B_{2L\epsilon}\setminus B_{L\epsilon}(p)\\[1.7ex]
              G_1+4\log{L\epsilon}-2\log{(1+\pi L^2)}-A_1(p_1)&others,
      \end{array}\right.$$
and
$$\phi_2=\left\{
     \begin{array}{ll}
       -\frac{w(\frac{x}{\epsilon})+2\log{(1+\pi L^2)}}{2}+2\log{L\epsilon}+
       \lambda_2r\cos\theta+\mu_2r\sin\theta
       +A_2(p)&x\in
        B_{L\epsilon}(p)\\[1.7ex]
       G_2-\eta H_2&x\in
        B_{2L\epsilon}\setminus B_{L\epsilon}(p)\\[1.7ex]
              G_2&others.
      \end{array}\right.$$
Here,
$$H_k=G_k-a_k\log r-A_k+\lambda_kr\cos\theta+\mu_kr\sin\theta,$$
and $\eta_i$ is a cut-off function which equals $1$ in
$B_{L\epsilon}(p)$, equals $0$ in $B_{2L\epsilon}^c(p)$.

Let $\Omega=\Sigma\setminus B_{L\epsilon}(p)$. By an argument
similar to the one used in Section 5, we can derive that
\begin{eqnarray*}
  \dint_\Sigma|\nabla\phi_1|^2dV_g&=&\dint_{B_{L\epsilon}(p)}|\nabla\phi_1|^2dxdy
      +\dint_\Omega|\nabla G_1|^2dV_g\\&&-2\dint_{\Sigma}\nabla G_1\nabla\eta H_1dV_g
      +\dint_{\Sigma}|\nabla\eta H_1|^2dV_g\\
   &=&\dint_{B_L}|\nabla w|^2dxdy-16\pi(L\epsilon)^2+\pi(L\epsilon)^2(\lambda_1^2+\mu_1^2)+
   O((L\epsilon)^4),
\end{eqnarray*}
\begin{eqnarray*}
  \dint_\Sigma|\nabla\phi_2|^2dV_g&=&\dint_{B_{L\epsilon}(p)}|\nabla\phi_1|^2dxdy
      +\dint_\Omega|\nabla G_2|^2dV_g\\&&-2\dint_{\Sigma}\nabla G_2\nabla\eta_2 H_2dV_g
      +\dint_{\Sigma}|\nabla\eta H_2|^2dV_g\\
   &=&\frac{1}{4}\dint_{B_L}|\nabla w|^2dxdy+8\pi(L\epsilon)^2+\pi(L\epsilon)^2(
    \lambda_2^2+\mu_2^2)+O((L\epsilon)^4),
\end{eqnarray*}
and
$$\begin{array}{ll}
   \dint_\Sigma\nabla\phi_1\nabla\phi_2dV_g&=
     \dint_\Omega\nabla G_1\nabla G_2dV_g+\dint_{B_{L\epsilon}(p)}
       \nabla\phi_1\nabla\phi_2dV_g\\[1.7ex]
     &\s-\dint\nabla G_1\nabla\eta H_2dV_g-
      \dint \nabla G_2\nabla\eta H_1dV_g+
      \dint_\Sigma\nabla\eta H_1\nabla\eta H_2dV_g\\[1.7ex]
    &=-\frac{1}{2}\dint_{B_L}|\nabla w|^2dxdy+\pi(L\epsilon)^2(\lambda_1\lambda_2
+\mu_1\mu_2)\\[1.7ex]
    &\s-4\pi(L\epsilon)^2+\dint_\Omega\nabla G_1\nabla G_2dV_g+O((L\epsilon)^4).
\end{array}$$
Note that
\begin{eqnarray*}
 \lefteqn{\int_\Omega(|\nabla G_1|^2+|\nabla G_2|^2+\nabla G_1\nabla
 G_2)dV_g}\\
        &=&\int_\Omega(|\nabla G_1|^2+|\nabla G_2|^2+\frac{\nabla G_1\nabla G_2+
          \nabla G_2\nabla G_1}{2})dV_g\\
    &=&-\int_{\partial B_{L\epsilon}(p)}
        (G_1\frac{\partial G_1}{\partial n}+G_2\frac{\partial G_2}{\partial n}
         +\frac{G_1\frac{\partial G_2}{\partial n}+
          G_2\frac{\partial G_1}{\partial n}}{2})dS_g\\&&
    +6\pi\dint_{B_{L\epsilon(p)}}(G_1+G_2)-6\pi\dint_\Sigma G_2dV_g.
\end{eqnarray*}
Applying Lemma 5.2, we get
$$\begin{array}{ll}
\dint_{\partial B_r(p)}G_1\frac{\partial G_1}{\partial n}dS_g=&
32\pi \log{r}-8\pi^2r^2+\pi r^2(
\lambda^2_1+\mu_1^2)\\[1.7ex]
& -8\pi A_1(p)+ 4\pi^2r^2A_1(p)-16\pi^2r^2\log{r}+O(r^4\log r),
\end{array}$$
$$\begin{array}{ll}
\dint_{\partial B_r(p)}G_2\frac{\partial G_2}{\partial n}dS_g=&
8\pi\log{r}+4\pi^2r^2+\pi r^2(
\lambda_2^2+\mu_2^2)\\[1.7ex]
& +4\pi A_2(p)+ 4\pi^2r^2A_2(p)+8\pi^2r^2\log{r}+O(r^4\log r),
\end{array}$$
$$\begin{array}{ll}
\dint_{\partial B_r(p)}G_1\frac{\partial G_2}{\partial n}dS_g=&
-16\pi \log{r}-8\pi^2r^2+\pi r^2(
\lambda_1\lambda_2+\mu_1\mu_2)\\[1.7ex]
& -8\pi A_2(p)+ 4\pi^2r^2A_2(p)+8\pi^2r^2\log{r}+O(r^4\log r),
\end{array}$$
$$\begin{array}{ll}
\dint_{\partial B_r(p)}G_2\frac{\partial G_1}{\partial n}dS_g=&
-16\pi \log{r}+4\pi^2r^2+\pi r^2(
\lambda_2\lambda_1+\mu_2\mu_1)\\[1.7ex]
& +4\pi A_1(p)+
4\pi^2r^2A_1(p)-16\pi^2r^2\log{r}+O(r^4\log r).
\end{array}$$
Note that
$$6\pi\dint_{B_{L\epsilon}(p)
 }(G_1+G_2)dV_g=6\pi^2 (L\epsilon)^2(A_1(p)+A_2(p))
-12\pi^2(L\epsilon)^2\log {L\epsilon}+6\pi^2 (L\epsilon)^2+O((L\epsilon)^4\log {L\epsilon}),$$
we get
\begin{eqnarray}\lefteqn{\dint_{\Sigma}(|\nabla\phi_1|^2+|\nabla\phi_2|^2
+\nabla\phi_1\nabla\phi_2)dV_g}\nonumber\\
&=&\frac{3}{4}\dint_{B_L}|\nabla w|^2dxdy-24\pi\log
{L\epsilon}\nonumber\\&&-6\pi A_1(p_1) -6\pi\dint
G_2dV_g+O((L\epsilon)^4\log{{L\epsilon}}).\end{eqnarray}
We  have
\begin{eqnarray*}
    \dint_\Sigma\phi_1dV_g&=&\epsilon^2\dint_{B_L}wdV_g-\dint_\Sigma\eta H_1dV_g
    +(4\log{{L\epsilon}}-2\log{(1+\pi L^2)}\\&&-A_1(p))(1-\dint_{B_{L\epsilon}(p)}dV_g)
    +\dint_{B_{L\epsilon}(p)}(\lambda_1x+\mu_1y)dV_g\\&&-\dint_{B_{L\epsilon}(p)}G_1dV_g
    +O((L\epsilon)^4\log {L\epsilon})\\[1.7ex]
    &=&\epsilon^2\dint_{B_L}we^{\varphi(\epsilon x,\epsilon y)}dxdy+4\log{{L\epsilon}}+2\pi (L\epsilon)^2
    \log{(1+\pi L^2)}
    \\&&-2\pi (L\epsilon)^2-A_1(p)-2\log(1+\pi L^2)+O((L\epsilon)^4\log{{L\epsilon}}).
\end{eqnarray*}
Because
$$\begin{array}{ll}
\dint_{B_{L\epsilon}(p)}G_2dV_g&=\dint_0^{L\epsilon}
 (2\log{r}+A_2(p))2\pi {L\epsilon}dr+O((L\epsilon)^4\log{L\epsilon})\\[1.7ex]
 &=2\pi r^2\log{r}-
   \pi r^2+(A_1(p)\pi r^2+O((L\epsilon)^4\log{L\epsilon}),
\end{array}$$
we can see that
$$\begin{array}{ll}
    \dint_\Sigma\phi_2dV_g&=\dint_\Sigma G_2dV_g-\dint_{B_{L\epsilon}(p)}G_2dV_g-
     \dint_\Sigma\eta\beta_2dV_g
     +(2\log{{L\epsilon}}+A_1(p))\dint_{B_{L\epsilon}}dV_g\\[1.7ex]
    &\s
    -\epsilon^2\dint_{B_L}\frac{w+2\log{(1+\pi L^2)}}{2}e^{\varphi(\epsilon x,
    \epsilon y)}dxdy+O((L\epsilon)^4\log {L\epsilon})\\[1.7ex]
    &=-\frac{\epsilon^2}{2}\dint_{B_L}we^{\varphi(\epsilon x,\epsilon y)}dxdy-\pi (L\epsilon)^2
    \log{(1+\pi L^2)}+\pi (L\epsilon)^2
    +O((L\epsilon)^4\log {L\epsilon}).
\end{array}$$
Hence,
\begin{equation}\dint(\phi_1+\phi_2)dV_g=-A_1(p)+4\log{L\epsilon}-\epsilon^2\log(1+\pi
L^2)-2\log(1+\pi L^2) +O((L\epsilon)^4\log {L\epsilon}).\end{equation}

Denote $B(p)=\frac{(b_1+\lambda_1)^2
+(b_2+\lambda_2)^2}{4}$, and $M=\frac{-\frac{K(p)}{2}+B(p)}{\pi}$,
we have
\begin{eqnarray}\log\dint_\Sigma e^{\phi_1}dV_g&=&\log{\epsilon^2}
+\epsilon^2M\log(1+\pi L^2)-2\epsilon^2(M+1)\log{L\epsilon}\nonumber\\&&+ O(\epsilon^2)+O((L\epsilon)^3\log
L)+O(\frac{1}{L^4}).\end{eqnarray} It is easy to see that
$$\dint_{B_{2L\epsilon}(0)}e^{\phi_2}dV_g=O((L\epsilon)^4),$$
and
$$\dint_{B_{2L\epsilon}(0)}e^{G_2}dV_g=O((L\epsilon)^4).$$
Since $\int e^{G_2}=1$, we get
\begin{equation}\log\dint_\Sigma e^{\phi_2}=\log(1-O((L\epsilon)^4))=O((L\epsilon)^4).
\end{equation}

In the end, we can deduce from (6.1), (6.2), (6.3) and (6.4) that
\begin{eqnarray*}
\Phi_0(\phi)&=&-4\pi-4\pi\log{\pi}+2\dint G_2dV_g-\epsilon^2(\log(1+\pi L^2)-2\log{L\epsilon})(1+
M)\\&&+O(\frac{1}{L^4})+O(\epsilon^3\log{L})+O((L\epsilon)^4\log{L\epsilon})+O(\epsilon^2).
\end{eqnarray*}
Let $L^4\epsilon^2=\frac{1}{\log(-\log\epsilon)}$. Then for
$\epsilon$ sufficiently small, we have
$$\Phi_0(\phi)<-4\pi-4\pi\log{\pi}+2\dint G_2dV_g.$$ This proves
our claim.

Therefore, if $\Sigma$ satisfies the condition that,
$$\max_{p\in\Sigma}K(p)<2\pi ,$$
 we can see that $u^\epsilon$ converges to $u^0=(u_1^0,u_2^0)$ in
$H_2:=H^{1,2}(\Sigma)\times H^{1,2}(\Sigma)$, hence it is clear
that $\Phi(u^0)=\inf_{u\in H_2}\Phi(u)$, that is, $u^0$ is a
minimizer of $\Phi_0=\Phi$. This completes the proof of the main
theorem.

\end{document}